\documentclass[axioms,communication,accept,moreauthors,12pt,a4paper]{mdpi} 

\lastpage{x}
\doinum{10.3390/------}
\pubvolume{xx}
\pubyear{2013}
\history{Received: xx / Accepted: xx / \\
Published: xx}

\usepackage{hyphenat}
\usepackage{enumitem}
\setenumerate{parsep=0pt,itemsep=0pt,leftmargin=.4in} 
\setitemize{parsep=0pt,itemsep=0pt,leftmargin=.4in}   
\setdescription{parsep=0pt,itemsep=0pt,leftmargin=.4in}

\newtheorem{de}{Definition}[section]

\newtheorem{re}[de]{Remark}

\newtheorem{te}[de]{Theorem}

\usepackage{soul}
\usepackage{amsfonts} 
\usepackage{amsmath}
\usepackage[all]{xy}

\newcommand{\ot}{\otg{B}}

\newcommand{\mproof}{\noindent{\bf Proof.}}
\newcommand{\twosid}[3]{\ar@<0.25ex>@{<-}[#1]^{#2}  \ar@<-1ex>[#1]_{#3}}

\def\ot{\otimes}

\Title{On transcendental numbers: new results and a little history} 

\Author{Solomon Marcus$^{1}$ and Florin F. Nichita$^{1,}$*}

\address{$^{1}$ Institute of Mathematics
of the Romanian Academy,
21 Calea Grivitei Street, 010702 Bucharest, Romania}

\corres{E-Mail: Florin.Nichita@imar.ro; \linebreak Tel:~(+40) (0) 21 319 65 06; Fax:~(+40) (0) 21 319 65 05 }

\abstract{Attempting to create a general framework for studying
new results on transcendental numbers, this paper 
begins with
a survey on 
transcendental numbers
and transcendence, it then presents several properties of the transcendental numbers $e$ and $\pi$,
and then it gives the proofs of new inequalities and identities for transcendental
numbers. 
Also, in relationship with these topics,
we study some implications for the theory of
the Yang-Baxter
equations, and we propose some open problems.}

\keyword{Euler's relation, transcendental numbers; transcendental operations / functions; transcendence; Yang-Baxter equation}

\MSC{ 01A05; 11D09; 11T23; 33B10; 16T25}

\begin{document}

\section{Introduction}

One of the most famous formulas in mathematics, the Euler's relation:
\begin{equation} \label{elapi2}
  e^{ \pi i} + 1 = 0 \ ,
\end{equation} 
contains the transcendental numbers $e$ and $ \pi$, the imaginary number $i$, the constants 0 and 1, and (transcendental) operations.
Beautiful, powerful and surprising, it has changed the mathematics forever.

We will unveil profound aspects related to it, and we will propose a counterpart:
\begin{equation} \label{new3}
 \vert e^i - \pi \vert <  e \ , 
\end{equation}
which has an interesting geometrical interpretation.

\begin{figure}[h!]
  \caption{{An interpretation of the inequality $ \  \vert e^i - \pi \vert <  e \ . $}}
  \centering
\setlength{\unitlength}{0.8cm}
\begin{picture}(4,3.2)
\thicklines
\put(1,0.5){\line(2,1){3}}
\put(4,2){\line(-2,1){2}}
\put(2,3){\line(-2,-5){1}}
\put(0.6,0.3){$A$}
\put(4.05,1.9){$B$}
\put(1.6,2.97){$C$}
\put(3.1, 2.7){${\bf 1}$}
\put(1, 1.7){${\bf e}$}
\put(2.8, 1.05){${\bf \pi}$}
\put(5,0.4){$\displaystyle
\widehat{B} \approx 1 \ radian$}
\end{picture}

\end{figure}
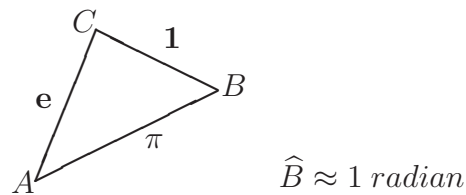

\bigskip




In some sense, the formulas (\ref{elapi2}) and (\ref{new3}) can be unified, and we will explain this in Section 3.

\bigskip

The next section of the current paper is a survey on
transcendental numbers, transcendental functions, transcendental operations and
transcendence (following the ideas from \cite{sm2, sm3}).
Section 3 deals with identities and inequalities of transcendental numbers.
In Section 4 we start with some comments,  we propose some open problems, and then we prove our new results.
Section 5 contains results related to the Yang-Baxter equation. 
In the last section we draw the final conclusions.

\section{Transcendence and transcendental numbers in mathematics}


The term ``transcendence'' has been introduced to express issues related to the Divinity. Implying initially the aforementioned attributes, it has produces
an inflation of uses nowadays: 
movies, novels, jazz albums, corporations, game brands, poems, internet websites, internet blogs, names of all kinds of doctrines, philosophies, 
etc.
The ``trans'' prefix comes from Latin, implying passing through a certain area, and it
is all around: translation, transport, transmission, transformation, transplant, transparent, transgress, transdisciplinarity etc.

\bigskip

One talked for the first time about {\em transcendence in mathematics} in the 18th century, and it was Leonhard Euler who initiated this discussion (see \cite{e}).
Leibniz was the first referring to transcendence in mathematics (see \cite{hb}).
Euler does not refer directly  to transcendental numbers (as it happened, subsequently, in the 19th century), but he refers to {\em transcendental operations}. 
The non-transcendental operations are: the addition, the subtraction, the multiplication, the division, the exponentiation and the rooting (applied to integers and to those numbers that are obtained 
from them by such operations). Why all these operations are considered to be non-transcendent? They are considered to be algebraic operations, because all numbers generated in this way 
are roots of some algebraic equations; we refer to those equations expressed by polynomial structures, where all mentioned operations are applied for finite number of times 
(focus on finite) and, as soon as this is not the case, we enter in the field of transcendence. Therefore, we are suddenly suggested that, in mathematics, {\em the idea of transcendence is 
essentially related to the idea of infinity}. But not to any kind of infinity. For Euler, the entry in mathematical transcendence is in the differential and integral calculus. 
Why is Euler considering the exponentiation a to x, the logarithm, the trigonometric sinus-cosine functions 
as being all {\em transcendental functions} and  operations? Because they are immediately associated, in a way or another, to the mathematical idea of integral, they are immediately associated 
to expressions including an integral (for instance, the logarithm is associated to the integral of dx/x). 
For Euler, the function log(x) + 7 is transcendent, but log(7) + x is algebraic.


The passing from non-transcendent to transcendent is gradual. 
Such an assertion deviates us from mathematical rigor, but it consolidates our understanding. We have in mathematics the situation of irrational numbers, 
where the word irrational would suggest transcendence, but it is immediately completed by the word algebraic the negation of transcendent, such as the square root of 2. 
We may consider the algebraic irrationals as a passing bridge, therefore intermediary when passing from non-transcendent to transcendent. It is not rigorous, but it is suggestive. 
We must add another thing: we talk about rational numbers (such as 2/ 3 or -3/7) and irrational numbers, and 
we are tempted to associate such irrational numbers as square root of 2 with something beyond reason. But the Latin etymology clarifies the fact that rational is sending here to ratio as fraction;
 so, in modern understanding, as a number of the form a/b, where a and b are integers, with b different from zero.


Joseph Liouville is the first who managed to encounter an example of {\em transcendental numbers}. 
This important discovery happened in the year 1844. Liouville introduced a class of real numbers wearing subsequently his name. 
A real number x is a Liouville number if there exist an integer b higher or equal to 2 and an infinite sequence of integers  ($a_1$, $a_2$, .., $a_n$, ...) so that x is 
the sum of the series having as general term the ratio between $a_k$ and the power of exponent k! of b. (For b = 10 and $a_k$ = 1 for any k, x becomes Liouville's constant.) 
Concerning his numbers, Liouville proves that they are not algebraic; 
they are transcendental. It was the first example of non-algebraic real numbers. Another presentation of {\em Liouville's numbers} stands on the manner of approximation by rational 
numbers: x is a Liouville number if for any natural number n there are integers p and q with q higher or equal to 2, such that the absolute value of the difference between x and p/q is strictly 
between zero and the fraction having 1 at the nominator and the value of q power n, at the denominator. It is thus clear that Liouville's numbers have the privilege of a tight approximation by rational numbers. 
This fact is against our intuitive expectations, because it shows that {\em in some respects transcendental numbers are nearer to rational numbers than algebraic irrationals}.


 The manner how the rational numbers are diligently „running” to get close to the transcendent numbers urges us to see the { \em transcendence process} 
in exactly the infinite range of approximations of transcendental numbers. This situation suggests us to generally regard {\em transcendence as result of an asymptotic process} 
when all stages are in the terrestrial universe, but the infinity of the number of stages makes impossible to be crossed in real time. 
Such a view increases the rightfulness of the hypothesis of absence of a sharp border between transcendent and terrestrial. 
 The relation between Liouville numbers and other important classes of real numbers has been recently studied by \cite{cs}. 
It is noticed that the Liouville numbers are encountered in remarkable classes of real numbers, but as rare exceptional phenomenon.


 There are other philosophical enigmas related to mathematical transcendence.  One would expect that transcendental numbers, such as $\pi$ or $e$, to be somehow farther from rational numbers, 
than irrational algebraic numbers. But in terms of approximation by rational numbers, it seems that this does not happen.
 It must be noticed, as well, that transcendental numbers are more than algebraic numbers, but „more” means here a metaphoric extension of such qualification, from finite to infinite. 
The mathematicians are expressing this by the following statement: ``The algebraic numbers form a countable set, while the transcendental numbers form an uncountable set; it is a set of the power of the continuum''.


\section{Transcendental numbers: identities and inequalities}

The following identities which contain the transcendental
numbers $ \ e $  and $ \ \pi $ are well-known:

\begin{equation} \label{pi3}
\int^{+ \infty}_{- \infty} e^{- x^2} dx = \sqrt{\pi} \ ,
\end{equation} 

\begin{equation} \label{pi4}
\int^{ + \infty}_{- \infty} e^{- i x^2} dx = \sqrt{\frac{\pi}{2}}(1-i) \ .
\end{equation} 

These formulas,  unified in the paper \cite{nams}, can be proved by contour integration.

\bigskip

Other inequalities for $ \ e $  and $ \ \pi $  (from \cite{ffn}) are quite new; for example, we list just two of them:
\begin{equation} \label{new6}
\mid e^{1-z} + e^{ \bar{z}} \mid > \pi \ \ \  \ \ \ \forall z \in \mathbb{C} \ , \  \ \ \ \ \ \ 
\end{equation}
$$  \int^{b}_{a} e^{- x^2} dx \  < \  \frac{e^e}{\pi} ( \frac{1}{e^{\pi a}}  -  \frac{1}{e^{\pi b}}) \ . $$

\bigskip

Still other relations will be proved in this paper: 
(\ref{new3}) and the experimental / informal formula (\ref{new4}).
\begin{equation} \label{new4}
 \vert \pi^i - i^{\pi} \vert = 2  \sin(\frac{\pi^2}{4} - \ln \sqrt{\pi})  \ .
\end{equation}


\bigskip

 Let us consider the two variable complex function
$ f : \mathbb{C} \times \mathbb{C} \rightarrow \mathbb{R} , \ \  f(z, w) = \vert e^z + e^w \vert $. 
The formulas (\ref{elapi2}), (\ref{new3}) and (\ref{new6}) can be unified using the function $f(z,w)$:\\
$ \ \ \ \ \ \ \ \ \ \ \ \ \ \ \ \  f(\pi i, \ 0) = 0 \ ,    \ \ \ \ \ \ \ \ \   \ \ \ \ \ \ \ \ \ \ \ \ f(1-z, \  \bar{z}) > \pi \ \  \forall z \in \mathbb{C} \ , \ \ \ \ \ \   \ \ \ \ \ \ \   \ \ \ \ \ \ \ \ \ \ \ f(i, \  \pi i + \ln \pi ) < e \ . $

\bigskip

\section{Further comments and proofs for our new results}

\bigskip

There exist real solutions for the equations
$$  x^2 - \pi x + (1+ \frac{1}{r})^r = 0\ , \ \ \ r \in \mathbb{Q^*} , $$
for r sufficiently small,
but
there are no real solutions for the ``limit'' equation
$$  x^2 - \pi x + e = 0, $$ 
because
$ \Delta= \pi^2-4e<0 \ .$
The question if $ \Delta= \pi^2-4e $ is a transcendental number is an open problem!

\bigskip

Resembling  the problem of squaring the circle,
the geometrical interpretation of the formula 
$\ \pi^2  < \ 4  e  \ $
could be stated as:
``The length of the circle with diameter $\pi$ is almost equal (and less) to the perimeter of a square with edges of length $e$''
 (see the Figure 2).
In this case, the area of the above circle is greater than the area of the above square, because
$ \pi^3 > 4 e^2 $.

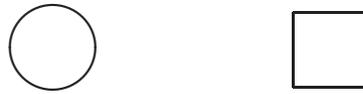
\begin{figure}[h!]
  \caption{{A circle with diameter $\pi$  and a square with edges of length $e$
have almost the same perimeter.}}
  \centering
\begin{picture}(220,70)
\thicklines
\put(100,40){\circle{31.4}}
\put(190,25){\line(0,1){28}}
\put(190,25){\line(1,0){28}}
\put(190,53){\line(1,0){28}}
\put(218,25){\line(0,1){28}}
\end{picture}
\end{figure}

{ \it
OPEN PROBLEMS (related to {\bf $\pi$}). For an arbitrary closed curve, we consider the smallest diameter ($d$)
and the maximum diameter (D). (These can be found by considering the center of mass of a body which
corresponds to the domain inside the
given curve.) 

(i) If $L$ is the length of the given curve
and the domain inside the
given curve is a convex set, then
we conjecture that:
$$ \frac{L}{D} \leq \pi \leq \frac{L}{d} \ \ .$$

(ii) Moreover, the first inequality becomes equality
 if and only if the second inequality becomes equality
if and only if the given curve is a circle.

(iii) If the area of the domain inside the
given curve is $A$, then $  \ \ \ d \ D  \ > A \ $.

(iv) The equation $ x^2 - \frac{L}{2} x + A =0 $ is not completely solved; for example,
if the given curve is an ellipse, solving this equation is an unsolved problem. }

\bigskip

\bigskip

Solving the equation 
$ \ \ \ $   $  \  x^2 \  + \  \sum_{0}^{\infty} 10^{-k!} \ x \  + \ \sum_{1}^{\infty} 10^{-k!}  \ = \ 0$,
in which two coefficients are
Liouville numbers, 
by using the
quadratic formula is hard;
however, one can
use the formula $ A x^2 + (A+C)x+C = (Ax+C)(x+1)$ in order to solve this equation.
 The solution $x=-1$ could be observed directly.

\bigskip
\bigskip

\begin{re}
{ \it The following informal proof for  the identity (\ref{new4}) 
is related to the
 computers working over the complex numbers.}

$ \ \ \  \vert \pi^i - i^{\pi} \vert = 
\vert e^{ i \ln \pi} - e^ {\pi \ln i} \vert =
\vert e^{ i \ln \pi} - e^ {i \pi \frac{\pi}{2}} \vert =
\sqrt{2 - 2( \cos (\ln \pi ) \cos (\frac{\pi^2}{2})
+ \sin (\ln \pi ) \sin (\frac{\pi^2}{2}))} =
\sqrt{2 - 2 \cos (\ln \pi - \frac{\pi^2}{2})} =
\sqrt{4 \sin^2 (\frac{\ln \pi - \frac{\pi^2}{2}}{2})} =
2  \sin(\frac{\pi^2}{4} - \ln \sqrt{\pi})  \  \approx 0.95  .$

\bigskip

In connection with the formula (\ref{new4}), one  could consider
the equation 
$$ x^i = i^x \ \ \ \ \ x \in \mathbb{R^*_+} \ , $$
which  is equivalent to
 $$  \ \ \ e^{\frac{\pi}{2}} = x^{\frac{1}{x}} \ \ \ \ \ x \in \mathbb{R^*_+} \ ,$$ 
and it has no real solution, because $ \frac{\pi}{2} > \frac{1}{e} $.

\end{re}

\bigskip

\bigskip

The proof for the inequality (\ref{new3}) follows in a similar fashion, using the Taylor series for $\cos x$
in order to approximate $\cos 1$,
and three digit approximations
for $e$ and $\pi$ (This proof is due to Dr. Cezar Joita).

\bigskip


\section{Transcendental numbers in mathematical physics}

In some special issues on
Hopf algebras, quantum groups and Yang-Baxter equations,
several papers \cite{a1, a2, a3, a4, a5, a6, a7, a8, a9, a10}, as well the feature paper
\cite{a11}, covered many topics related to the Yang-Baxter equation, ranging from mathematical physics
to Hopf algebras, and from Azumaya Monads to quantum computing.

The terminology of this section is compatible with the above cited papers.
 $V$ is a complex vector space, and $ \  I_j : V^{\ot j} \rightarrow 
 V^{\ot j} \ \ \ \forall j \in \{1, 2 \} $ are identity maps.
We consider 
 $ \  J : V^{\ot 2} \rightarrow 
 V^{\ot 2} $ a linear map which satisfies
$ J \circ J = - I_2 \ $ and $ \ \ \ J^{12} \circ J^{23} = J^{23}  \circ J^{12} \ , \  $ where
$ J^{12} = J \ot I_1 \ , \ \ J^{23} = I_1 \ot J$.

Then,  
$R (x) = \cos x I_2 + \sin x J$
satisfies the colored Yang-Baxter equation:
\begin{equation} \label{yb}
R^{12}(x) \circ  R^{23}(x+y)  \circ  R^{12}(y) \ =  R^{23}(y)  \circ  R^{12}(x+y)  \circ  R^{23}(x) \ .
\end{equation}



The proof of (\ref{yb}) could be done by writing
$R (x) = e^{x \ J}$, and checking that (\ref{yb}) reduces to

$ x \ J^{12} \ + \ (x+y) \ J^{23} \ + \ y J^{12} = 
y \ J^{23} \ + \ (x+y) \ J^{12} \ + \ x J^{23} $.


Such an operator $J$ could have, in dimension two, the following matrix form (for $ \alpha \in \mathbb{R}$):

\begin{equation} \label{rmatcon2}
\begin{pmatrix}
0 & 0 & 0 & \frac{1}{\alpha} i\\
0 & 0 & i & 0\\
0 & i & 0  & 0\\
\alpha i & 0 & 0 & 0
\end{pmatrix}
\end{equation}

Based on results from the previous section, a counterpart for the formula 
$$ \ \ e^{\pi \ J} + I_4 = 0_4 \ \ \  \ \ \ \  J, \ I_4 , 0_4 \in \mathcal{M}_4 ( \mathbb{C}) $$ 
could be the following inequality:
\begin{equation} \label{new7}
X^2 + e I_2 >  \pi X \ \ \ \ \ \ \ \ \ \ \ \ \ \ \ \ \ \ \ \ \ \ \forall \  X \in \mathcal{M}_2 ( \mathbb{R^*_+)})   ,  \ trace(X)> \pi \ .
\end{equation}

\bigskip

Replacing the above condition $J^{12} \circ J^{23} = J^{23}  \circ J^{12}$ with 
$J^{12} \circ J^{23} = - J^{23}  \circ J^{12}$, the authors of \cite{majo} obtained interesting results
 (a new realization of doubling degeneracy based on emergent Majorana operator, new solutions for the Yang-Baxter equation, etc).
For example, in dimension two, the matrix form of this new operator $ J $ could be:
\begin{equation} \label{rmatcon2}
\begin{pmatrix}
0 & 0 & 0 & 1\\
0 & 0 & 1 & 0\\
0 & -1 & 0  & 0\\
-1 & 0 & 0 & 0
\end{pmatrix}
\end{equation}

With this case we enter into the world of  the quaternions and Clifford algebras. 

\begin{te}
 For a Boolean algebra, the map $R(a, \ b) = (a \rightarrow b, \ a)$
is a solution for the constant Yang-Baxter equation:
\begin{equation} \label{cyb}
R^{12} \circ  R^{23}  \circ  R^{12} \ =  R^{23}  \circ  R^{12} 
 \circ  R^{23} \ .
\end{equation}
\end{te}

{ \bf Proof.} The proof is direct. 
\bigskip

\section{Conclusions}

The author of \cite{BirdsFrogs} considers two types of scientists:
{\bf birds} (they resemble
scientists with a broad vision, who try to unify theories, who obtain results of interest for a large readership) and
{\bf frogs} (which are less influential). 
Solomon Marcus 
 used the terms of Francis Bacon (Novum Organum),
 {\bf bees} versus {\bf ants},  
while describing
mathematicians who are involved in 
many different areas of research 
versus the mathematicians who work on a restricted domain.

Because there is a huge number 
of new disciplines, it is important to have a transdisciplinary understanding of the world:
a transdisciplinary approach (see \cite{b1, b2, fn}) 
attempts to discover what is between disciplines, across different 
disciplines, and beyond all disciplines. Transcendence is a concept 
which plays an important role in theology, in science and in art; it can be
considered beyond all disciplines. The Yang-Baxter equation appears across
different disciplines. Mathematical Physics is at the border of two disciplines.

\bigskip


Our paper is written in transdisciplinary fashion.
We used results and concepts
from algebra, mathematical analysis, mathematical physics, geometry, history of mathematics,
numerical analysis, epistemology, philosophy etc. 
Attempting to continue the approach of our recent papers and talks on the transcendental numbers (see \cite{sm2}, \cite{sm3}, 
\cite{ffn}),
we brought together our investigations and we presented new results.


\bigskip
\bigskip
\begin{center}
{\bf Acknowledgements}

We would like to thank Dr. Cezar Joita for his participation to our investigations.

\end{center}
\bigskip

\bibliographystyle{mdpi}

\begin{thebibliography}{----}

\bibitem{sm2} Marcus, S. {\em Transcendenta ca paradigma universala} (in Romanian), Convorbiri Literare {\bf 2014}, February, 17-28.

\bibitem{sm3} Marcus, S. {\em Transcendence, as a universal paradigm},
BALANCE, A Club of Rome Magazine {\bf 2015}, no.1, to appear.

\bibitem{e} Petrie, B.J. {\em Leonhard Euler’s use and understanding of mathematical transcendence}, Historia Mathematica {\bf 2012}, 39, 280-291.

\bibitem{hb} Breger, H. {\em Leibniz Einfuhrung des Transzendenten}, 300 Jahre „Nova Methodus” von G.W. Leibniz (1684-1984). Studia Leibnitiana {\bf 1986}, Sonderheft XIV, 119-132.



\bibitem{cs} C.S. Calude, C.S.; Staiger, L.
{\em Liouville numbers, Borel normality and algorithmic randomness}, CDMTCS-448, Auckland, December 2013.

\bibitem{nams} Desbrow, D. {\em On Evaluating {$\int^{+ \infty}_{- \infty} e^{ax(x-2b)} dx$} by Contour Integration Round a Parallelogram}, The Amer. Math. Month. {\bf 1998}, Vol 105, Number 8, 726-731.


\bibitem{ffn} Nichita, F.F. {\em On Transcendental Numbers}, Axioms {\bf 2014}, 3(1), 64-69.

\bibitem{a1} Kanakoglou, K. {\em Gradings, Braidings, Representations, Paraparticles: Some Open Problems}, Axioms {\bf 2012},{ 1(1)}, 74-98.

\bibitem{a3} Underwood R.G. {\em Quasi-triangular Structure of 
Myhill-Nerode Bialgebras}, Axioms. {\bf 2012}, { 1(2)}, 155-172.

\bibitem{a5} Schmidt, J.R. {\em From Coalgebra to Bialgebra for the Six-Vertex Model: The Star-Triangle Relation as a Necessary Condition for Commuting Transfer Matrices}. Axioms {\bf 2012}, { 1(2)}, 186-200.

\bibitem{a6} Nichita, F.F.; Zielinski, B. {\em The Duality between Corings and Ring Extensions}. Axioms {\bf 2012}, { 1(2)}, 173-185.

\bibitem{a2} Links, J. 
{ \em Hopf Algebra Symmetries of an Integrable Hamiltonian for Anionic Pairing}, Axioms {\bf 2012}, { 1(2)}, 226-237.

\bibitem{a4} Hoffnung, A.E. {\em The Hecke Bicategory}, Axioms {\bf 2012},{ 1(3)}, 291-323.

\bibitem{a7} Nichita, F.F. { \em
Yang-Baxter systems, algebra factorizations and
braided categories}, Axioms {\bf 2013}, { 2(3)}, 437-442.



\bibitem{a8} Lebed, V. {\em R-Matrices, Yetter-Drinfel$'$d Modules and
 Yang-Baxter Equation}, Axioms {\bf 2013}, { 2(3)}, 443-476.

\bibitem{a9}
Iordanescu, R.; Nichita, F.F.; Nichita I.M. {\em The Yang-Baxter Equation, (Quantum) Computers and Unifying Theories}, Axioms {\bf 2014}, 3(4), 360-368.

\bibitem{a10} Mesablishvili, B.; Wisbauer, R. {\em Azumaya Monads and Comonads}, Axioms {\bf 2015}, 4(1), 32-70.

\bibitem{a11} Nichita, F.F. { \em Introduction to the Yang-Baxter Equation with Open Problems}, Axioms {\bf 2012}, {  1(1)}, 33-37.

\bibitem{majo} Yu, L.W.; Ge, M.L. {\em  More about the doubling degeneracy operators associated with Majorana fermions and Yang-Baxter equation},
    Scientific Reports {\bf 2015}, 5,
Article number: 8102, 7 pages.

\bibitem{BirdsFrogs} Dyson, F. {\em Birds and Frogs}, Notices of the AMS {\bf 2009},
Vol. 56, No. 2, 212-223.


\bibitem{b1} Nicolescu, B. { \em 
Manifesto of Transdisciplinarity}, State University of New York (SUNY) Press, New York, translation  in English by Karen-Claire Voss, 2002.

\bibitem{b2} Nicolescu, B. { \em Transdisciplinarity - 
 past, present and future}, in Moving Worldviews - Reshaping sciences, policies and practices for endogenous sustainable development {\bf 2006}, COMPAS Editions, Holland, edited by Bertus Haverkort and 
Coen Reijntjes, 142-166.

\bibitem{fn} Nichita, F.F. {\em On Models for Transdisciplinarity}, Transdisciplinary Journal of 
Engineering and  Science {\bf 2011}, Vol. 2011, 42-46.




\end{thebibliography}
\makeatletter
\makeatother

\end{document}